\documentclass[a4paper,12pt,doublespace]{article}
\usepackage{epsfig,graphicx}
\usepackage{makecell}
\usepackage{amssymb,bbm}
\usepackage{amsmath}
\usepackage{amsfonts}
\usepackage{setspace}
\usepackage[toc,page]{appendix}
\usepackage{enumerate}
\usepackage{subfig}

\usepackage{float}
\usepackage[dotinlabels]{titletoc}
\newcounter{theorem}

\newtheorem{definition}{Definition}

\title{Some Examples of Graphs Suggesting That the Discrete Curvature Does Sense the Smooth One\date{}}
\author{G\"{o}k\c{c}e \c{C}AKMAK \footnote{Eskisehir Technical University, Science Faculty, Department of Mathematics, 26470, Eski\c{s}ehir, Turkey.  e-mail: gokcecakmak@eskisehir.edu.tr}
\thanks{Corresponding Author.}
\and Ali DEN\.{I}Z \footnote{Eskisehir Technical University, Science Faculty, Department of Mathematics, 26470, Eski\c{s}ehir, Turkey. e-mail: adeniz@eskisehir.edu.tr } \and  \c{S}ahin KO\c{C}AK \footnote{Anadolu University (emeritus), Science Faculty, Department of Mathematics, 26470, Eski\c{s}ehir, Turkey. e-mail: skocak@anadolu.edu.tr} \and Murat L\.IMONCU \footnote{Eskisehir Technical University, Science Faculty, Department of Mathematics, 26470, Eski\c{s}ehir, Turkey.  e-mail: mlimoncu@eskisehir.edu.tr}}
\begin{document}
\maketitle
\thispagestyle{empty}

\begin{abstract}
In this note, using some regular triangular tilings of the sphere, the Euclidean plane and the hyperbolic plane, we examine the potential relationship between their discrete Bakry–\'{E}mery curvatures and the smooth curvatures of their ambient space forms.
\end{abstract}


\section{Introduction}
The notion of graph curvature has attracted considerable interest in recent years (\cite{CushingCalc2022}, \cite{CushingEig2022}, \cite{Fathi2022}, \cite{Klartag2016}, \cite{Yau2011}, \cite{Viola2021}). There are several approaches to this discrete curvature notion and one of them, the so-called Bakry-\'{E}mery graph curvature, has direct roots in smooth differential geometry (\cite{Bauer2017}). We wondered whether the Bakry-\'{E}mery curvature of a graph somehow reflects the smooth curvature of an ambient surface in which the graph is embedded in such a way that the edges are realized as geodesic segments on the surface.

\begin{figure}[H]
  \centering
  \includegraphics[scale=0.7]{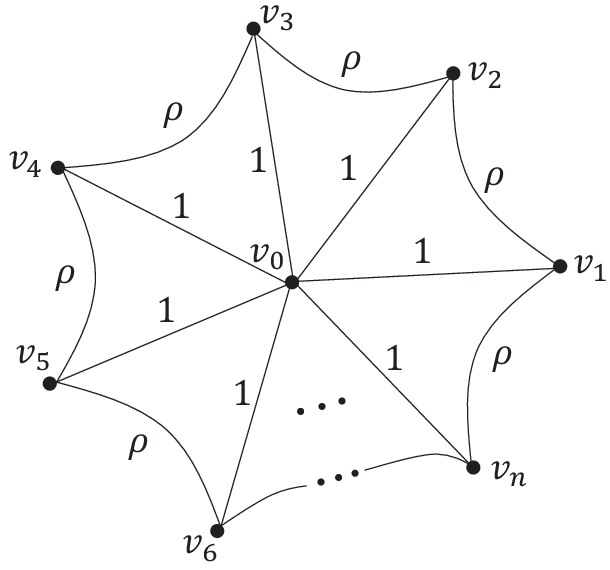}
   \caption{The umbrella graph $G_{n, \rho}$} \label{sekil1}
\end{figure}

In a previous attempt \cite{Cakmak2023}, we considered a family of weighted graphs $G_{n,\rho}$ (we called the umbrella graphs, see Figure \ref{sekil1}) depending on two parameters $n$ and $\rho$ with $n\geq 3$ and $0<\rho <2$. Given $n$ and $\rho$, there exists a unique space form (i.e. a sphere, the Euclidean plane or a hyperbolic plane) with a well-defined curvature such that the graph $G_{n,\rho}$ can be geodesically embedded into. Our hope was that the discrete Bakry–\'{E}mery curvature of the weighted graph $G_{n,\rho}$ (at the central vertex) would somehow reflect the curvature of the corresponding space form. To compute the discrete curvature, we used a weighted version of the discrete curvature setting (see \cite{CushingEig2022}) whereby as weights of the edges we took the edge-lengths and as vertex measures the sums of the lengths of edges emanating from vertices. We then compared the weighted discrete curvature of the graph $G_{n,\rho}$ with the smooth curvature of the corresponding ambient space form (whose curvature can be easily computed in dependence of $\rho$). The result was rather disappointing; the main issue was that we obtained positive values for the discrete curvature for all $G_{n,\rho}$, even for those whose ambient spaces were hyperbolic. Nevertheless, there was one interesting feature pointing still to a relationship between discrete and smooth curvatures: If we fix $n$ and increase $\rho$ in such a way that $G_{n,\rho}$ is embedded successively into the sphere, Euclidean plane and hyperbolic plane, the discrete curvature of $G_{n,\rho}$ clearly decreases, though remaining positive and different from the corresponding smooth curvature. Although this tendency could still be interpreted as a hint that the discrete curvature “tries” to sense the smooth curvature, it must be admitted that this state of affairs is unsatisfactory. There could be two reasons for this displeasing picture: Either the chosen procedure for the computation of the weighted discrete curvature is not quite appropriate or the used test-graphs $G_{n,\rho}$ are not well-chosen in the sense that any vertex is within 1-distance from the central vertex, but genuine 2-neighbourhoods are crucial for computation of discrete curvatures. We will address these two issues, which are technically rather involved, in another note.

In the present note, we want to use another family of graphs for testing a possible relationship between discrete and smooth curvatures. We will consider regular triangular tilings of the sphere, the Euclidean plane and the hyperbolic plane and choose the corresponding smooth curvatures of the space forms such that the edges of the triangulations are of unit length. Then, we will compute the unweighted discrete curvature of the graph in the sense of Klartag et. al. \cite{Klartag2016}.

We now fix the definitions and notations \cite{CushingEig2022, Klartag2016, Viola2021}:
Given a finite, simple, connected graph $G=(V, E)$, the Laplacian $\Delta$ acting on $f:V\to \mathbb{R}$ is defined by
\[ \Delta f(x)=\displaystyle\sum_{v\sim x}(f(v)-f(x)),  \]
where $v\sim x$ denotes that $v$ is a neighbour of $x$ in $G$.

The Laplacian gives rise to the symmetric bilinear forms
\begin{align*}
  2\Gamma(f,g) &:=\Delta(f,g)-f\Delta g-g\Delta f , \\
   2\Gamma_{2}(f,g)  & :=\Delta(\Gamma(f,g))-\Gamma(f, \Delta g)-\Gamma(g, \Delta f).
\end{align*}

$\Gamma(f,g)$ can be expressed in more explicit terms as
\[ \Gamma(f,g)(x)=\dfrac{1}{2}\displaystyle\sum_{v\sim x}(f(v)-f(x))(g(v)-g(x)). \]
We will use the abbreviations $\Gamma(f):=\Gamma(f,f)$ and $\Gamma_{2}(f):=\Gamma_{2}(f,f)$.

\begin{definition}
The Bakry-\'{E}mery curvature for the dimension $\infty$ at a vertex $x\in V$ of a graph $G=(V, E)$ is the maximum value $K\in \mathbb{R} \cup \{-\infty\}$ such that for any real function $f:V\to \mathbb{R}$, $\Gamma_{2}(f)(x)\geq K\Gamma(f)(x)$.
\end{definition}

We remark that the Bakry-\'{E}mery curvature $K_G(x)$ of the graph $G$ at the vertex $x$ can also be expressed as
\[
K_G(x)=\inf_{f} \dfrac{\Gamma_{2}(f)(x)}{\Gamma(f)(x)}, \qquad \qquad (\Gamma(f)(x)\neq 0).
\]
This can be shown along the lines of \cite{Viola2021} and it results from the fact that for $\Gamma(f)(x)=0$, the quantity $\Gamma_{2}(f)(x)$ is nonnegative. For computations, one can assume without loss of generality $f(x)=0$ and in the following computations we will use the program of Cushing et. al. \cite{CushingCalc2022}.

In the next section, we will specify the considered regular triangular tilings of the space forms and determine the curvature of each space form such that the edge-lengths of its triangular tiling are of unit length. As the edge-lengths of a graph are assumed to be unit in the setting of unweighted discrete graph curvature, we wanted to be consistent with that assumption. This convention uniquely determines for each tiling the curvature of the corresponding space form. We then compute the discrete curvatures of the graphs associated with these tilings (i.e. of the graphs of their one skeletons) and compare them with the smooth curvatures of the corresponding space forms.

\section{The Curvature of the Regular Triangular Tilings}
We describe and illustrate in the following the considered graphs (as one skeletons of certain tilings) and indicate how the smooth curvatures of their ambient space forms are computed.

\begin{figure}[H]
\begin{minipage}{.32\linewidth}
\centering
\subfloat[]{\label{figure1a}\includegraphics[scale=.35]{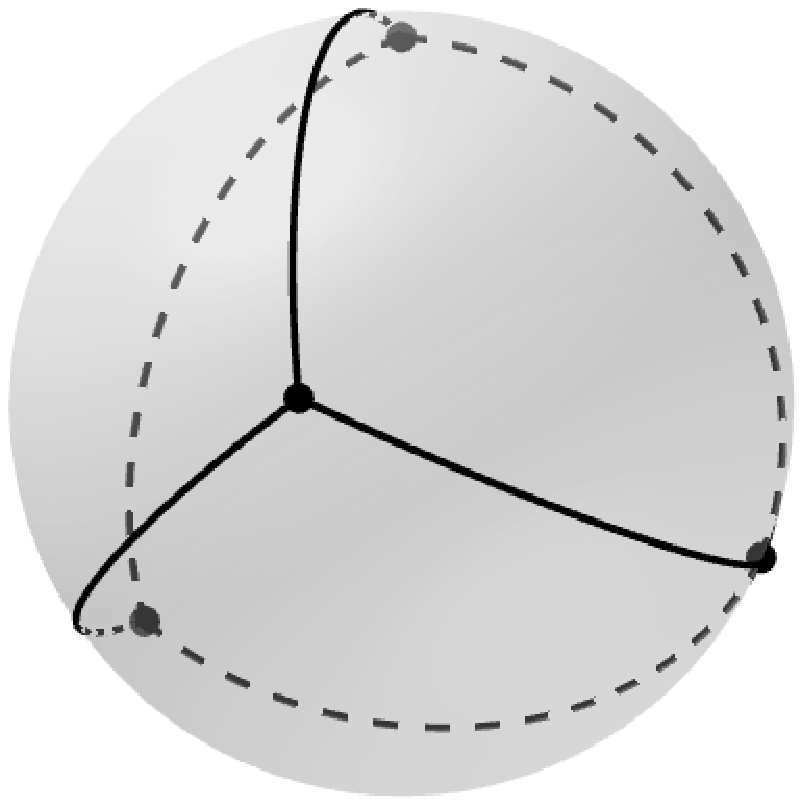}}
\end{minipage}%
\begin{minipage}{.35\linewidth}
\centering
\subfloat[]{\label{figure1b}\includegraphics[scale=.3]{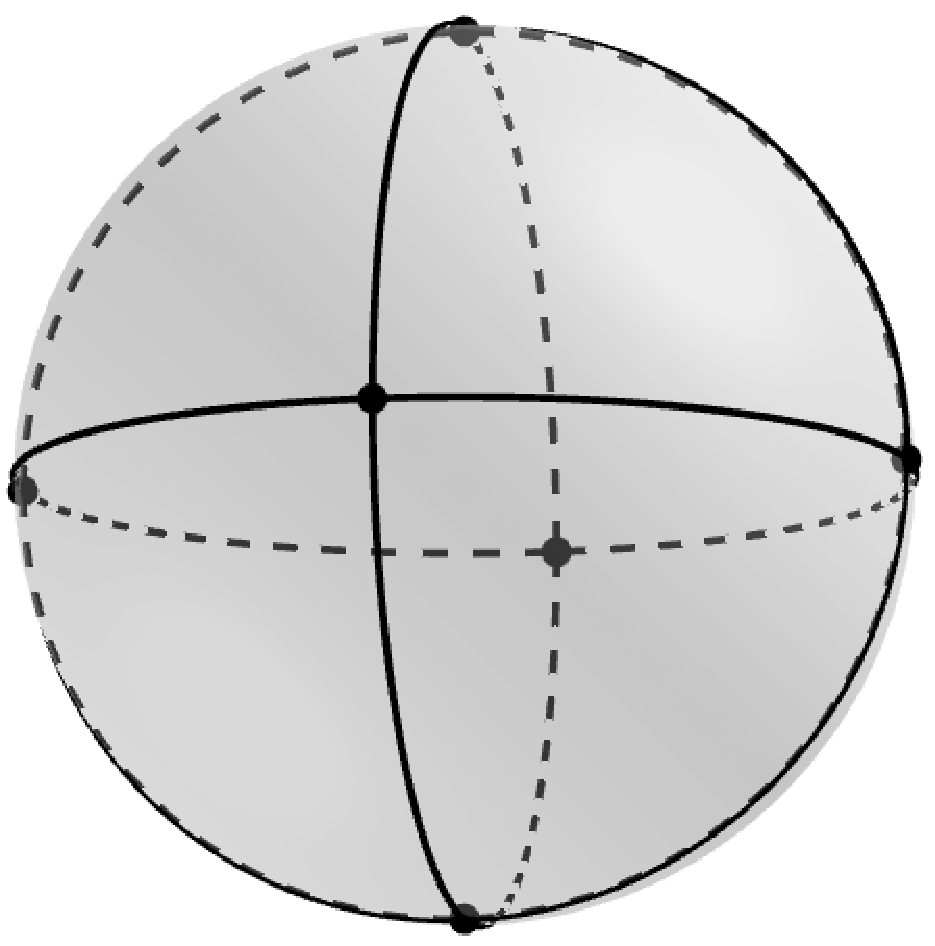}}
\end{minipage}%
\begin{minipage}{.32\linewidth}
\centering
\subfloat[]{\label{figure1c}\includegraphics[scale=.35]{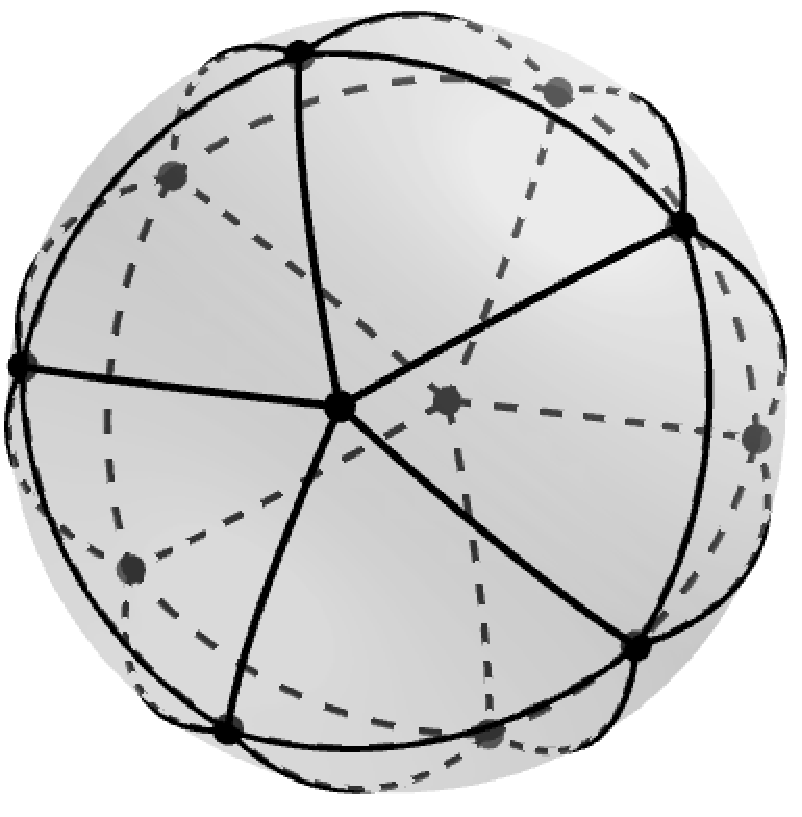}}
\end{minipage}\par\medskip
\begin{minipage}{.32\linewidth}
\centering
\subfloat[]{\label{figure1d}\includegraphics[scale=.1]{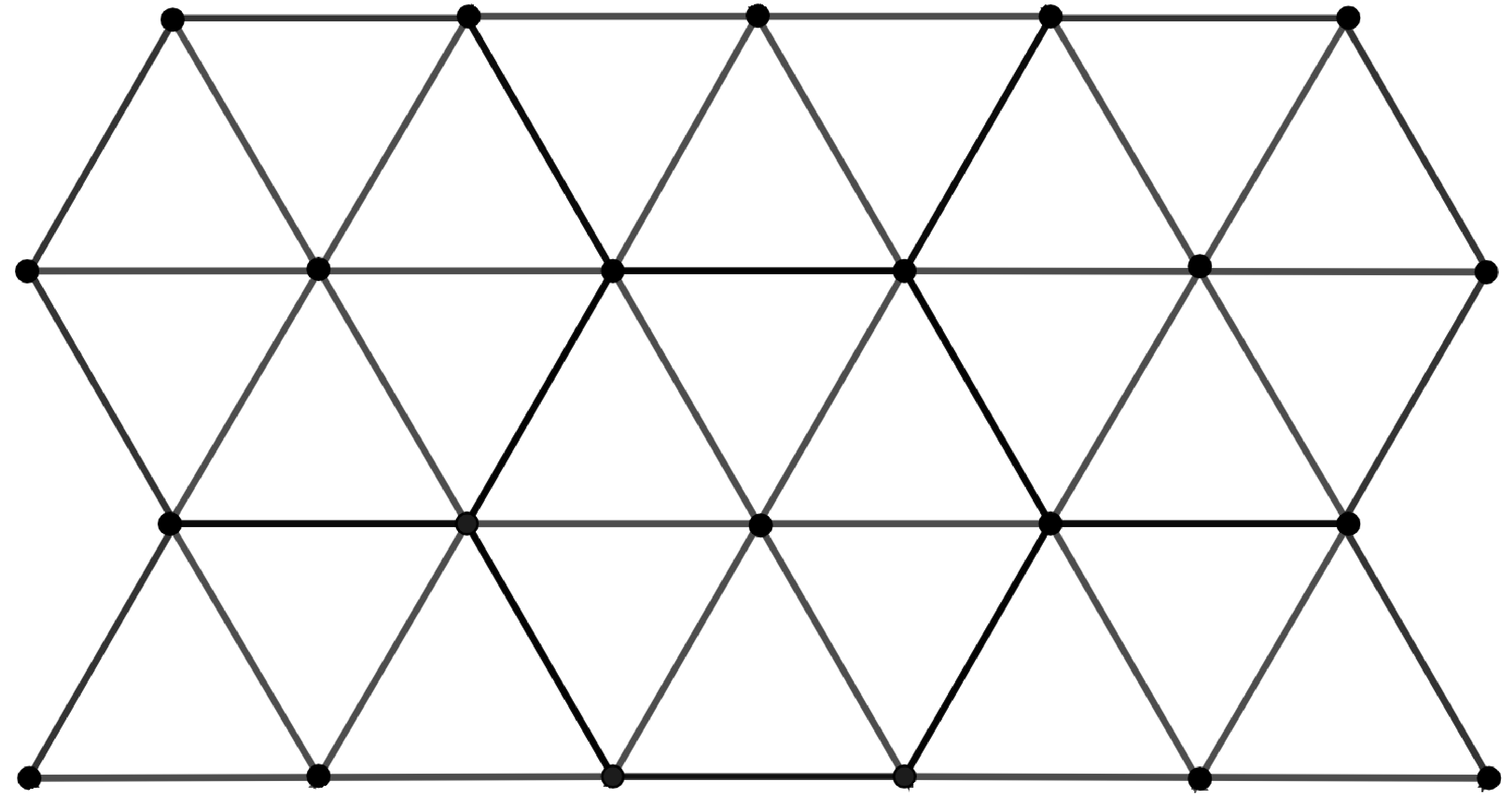}}
\end{minipage}%
\begin{minipage}{.35\linewidth}
\centering
\subfloat[]{\label{figure1e}\includegraphics[scale=.05]{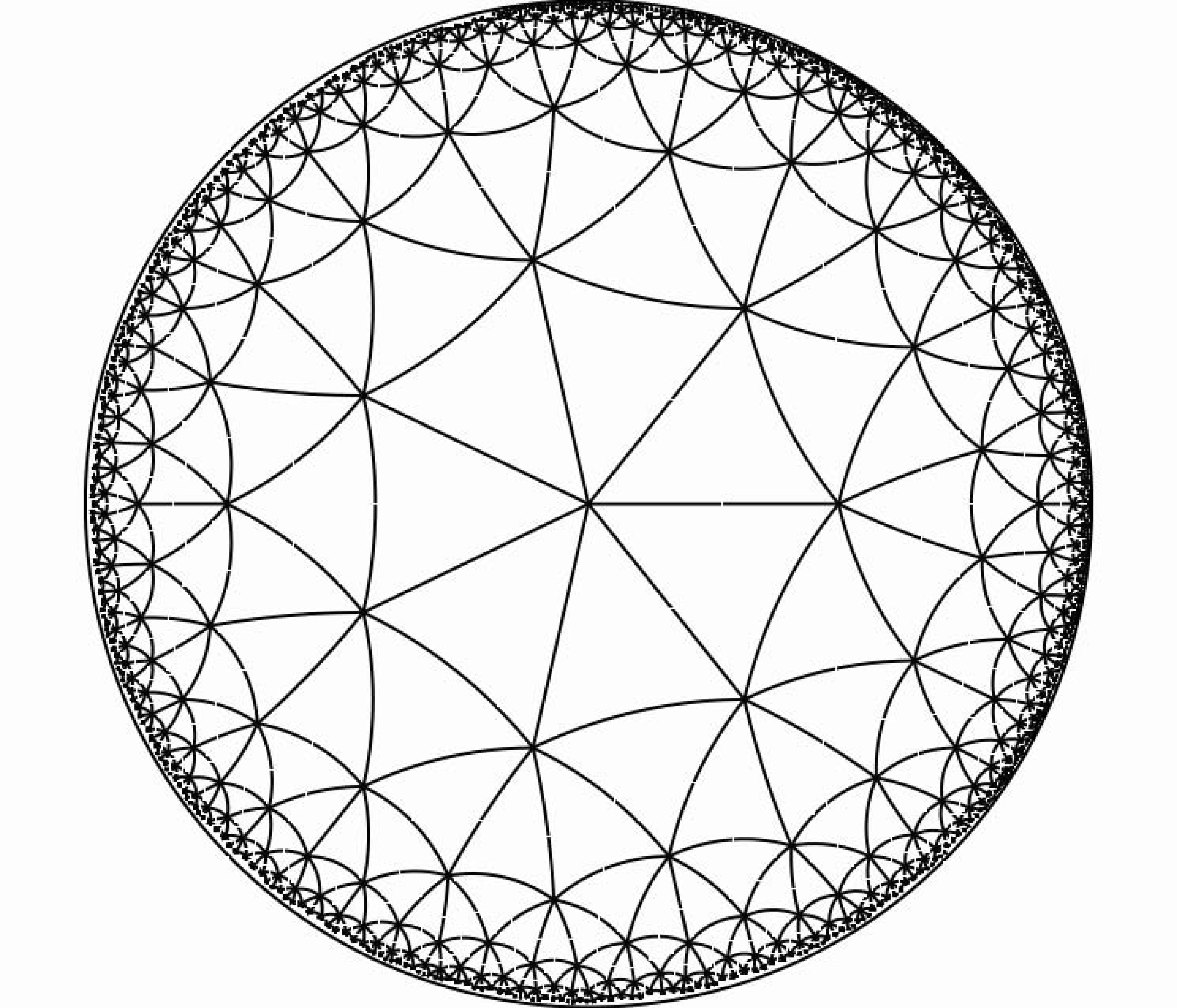}}
\end{minipage}%
\begin{minipage}{.32\linewidth}
\centering
\subfloat[]{\label{figure1f}\includegraphics[scale=.05]{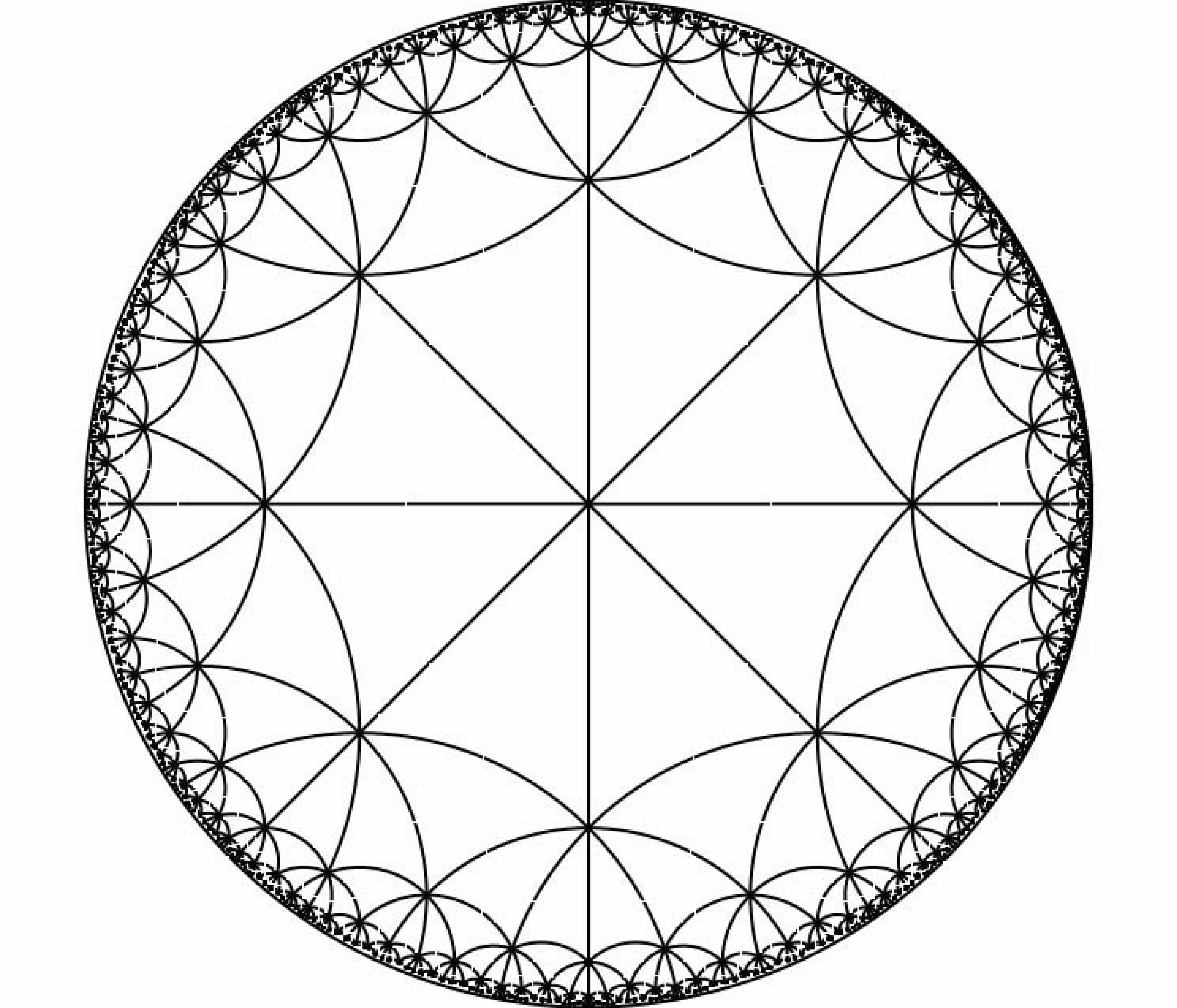}}
\end{minipage}\par\medskip
\centering
\subfloat[]{\label{figure1g}\includegraphics[scale=.05]{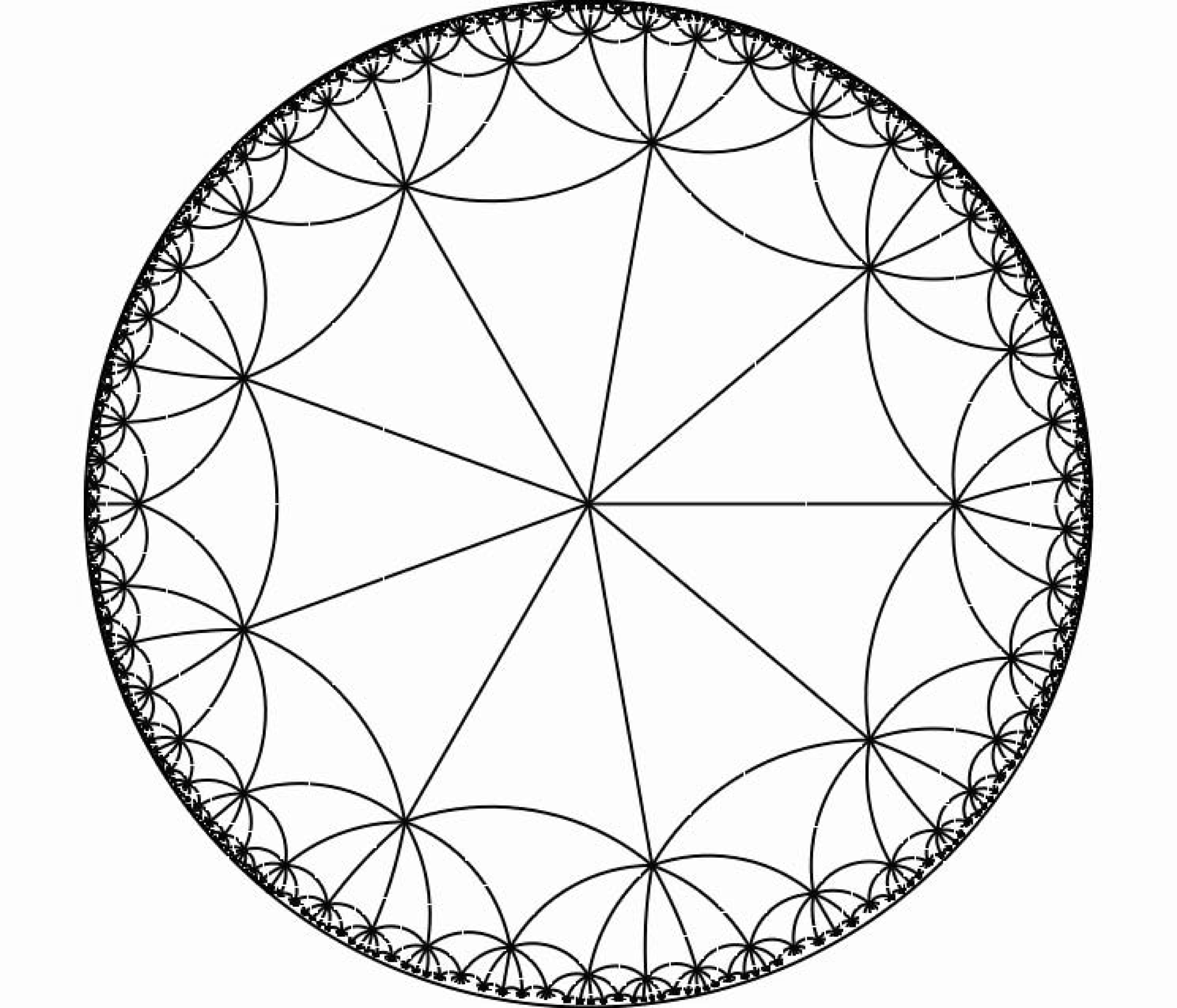}}
\caption{The order-3, -4, and -5 regular triangular tilings of the sphere are represented in (a), (b), (c), the order-6 regular triangular tiling of the plane is represented in (d) and the order-7, -8, and -9 regular triangular tilings of the hyperbolic plane are represented in (e), (f), (g).}
\label{figure1}
\end{figure}

In the case of spherical tilings we determine radii of the spheres so that the corresponding tilings have unit edge lengths. We use thereby the spherical formula $\cos \frac{1}{R}=\cos^2 \frac{1}{R}+\sin^2 \frac{1}{R} \cos \alpha$ for the unit-length equilateral triangle, where $\alpha$ is $2\pi/3$ for the order-3 tiling, $\alpha=\pi/2$ for the order-4 tiling and $\alpha=2\pi/5$ for the order-5 tiling (see Figure \ref{figure1a}, \ref{figure1b}, \ref{figure1c}). The formula gives $R=1/ \arccos (-1/3)\approx 0.523$ for $\alpha=2\pi/3$, $R=2/ \pi\approx 0.636$ for $\alpha=\pi/2$ and $R=1/ \arccos (1/\sqrt5)\approx 0.903$ for $\alpha=2\pi/5$. The corresponding curvatures of the spheres are $\kappa=1/R^2\approx 3.652$ for $\alpha=2\pi/3$, $\kappa\approx 2.467$ for $\alpha=\pi/2$ and $\kappa\approx 1.226$ for $\alpha=2\pi/5$.

In the case of hyperbolic tilings, we determine the curvatures $\kappa$ of the ambient hyperbolic planes by the formula $\cosh \sqrt{-\kappa}=\cosh^2 \sqrt{-\kappa}-\sinh^2 \sqrt{-\kappa} \cos\alpha$ for the unit-length equilateral triangle. From this formula, one can obtain the approximate values $\kappa\approx -1.189$ for the order-7 tiling (with $\alpha=2\pi/7$), $\kappa\approx -2.337$ for the order-8 tiling (with $\alpha=\pi/4$) and $\kappa\approx -3.441$ for the order-9 tiling (with $\alpha=2\pi/9$) (see Figure \ref{figure1e}, \ref{figure1f}, \ref{figure1g}).

For the Euclidean plane, we take simply the order-6 regular triangular tiling with unit edge length (see Figure \ref{figure1d}).

\begin{figure}[H]
\centering
\scalebox{.9}{
\begin{minipage}{.32\linewidth}
\centering
\subfloat[]{\label{figure2a}\includegraphics[scale=.1]{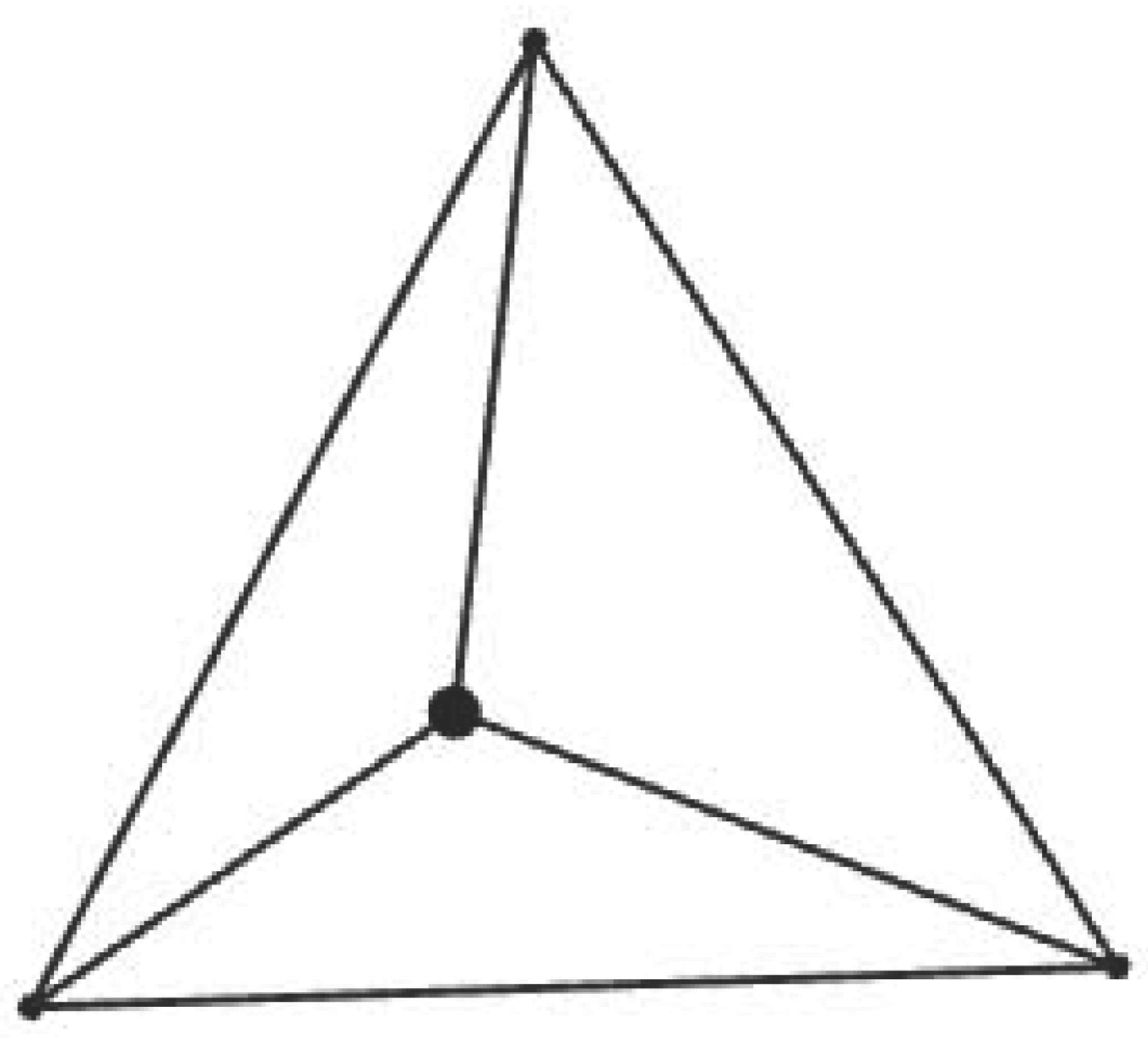}}
\end{minipage}%
\begin{minipage}{.35\linewidth}
\centering
\subfloat[]{\label{figure2b}\includegraphics[scale=.1]{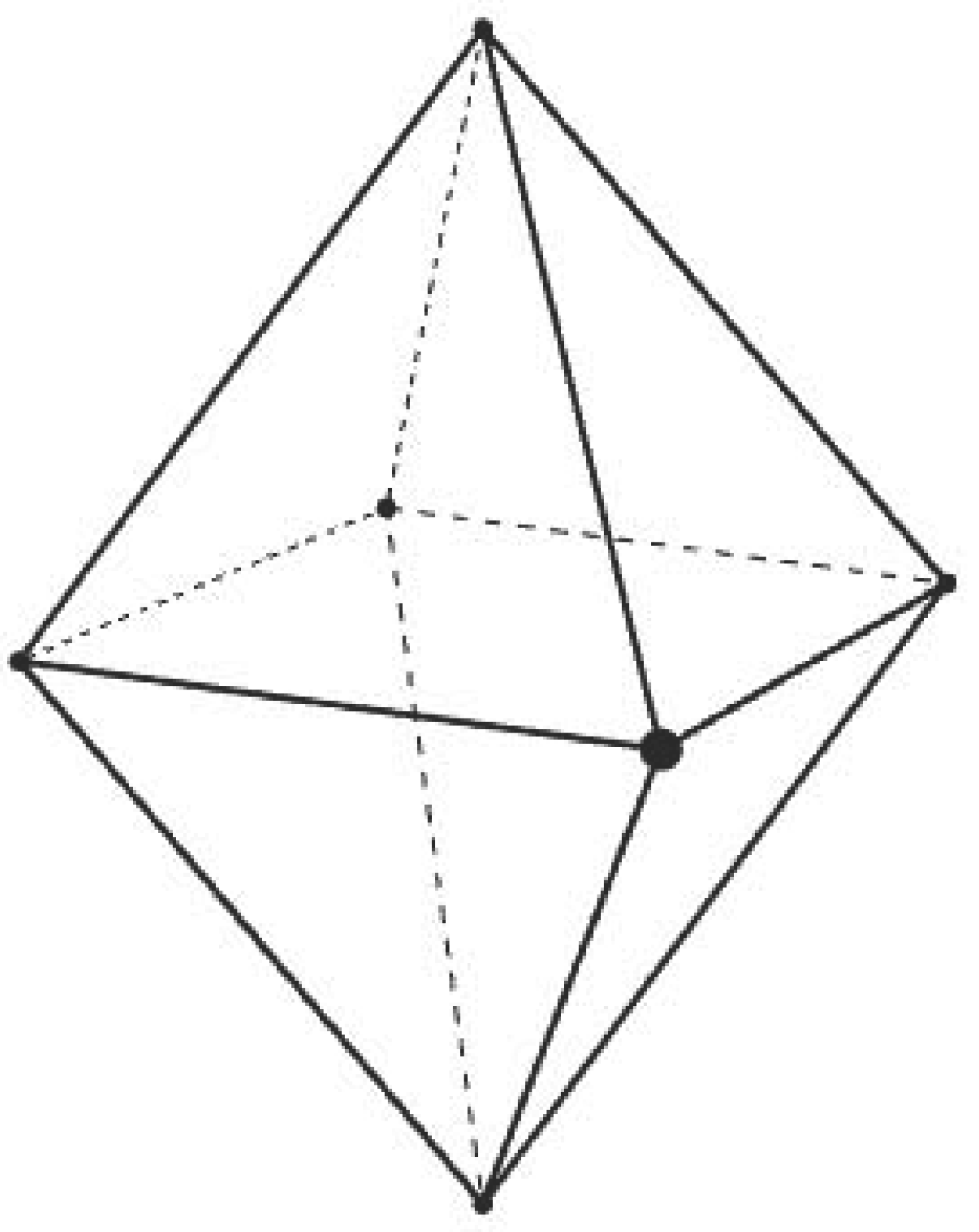}}
\end{minipage}
\begin{minipage}{.32\linewidth}
\centering
\subfloat[]{\label{figure2c}\includegraphics[scale=.1]{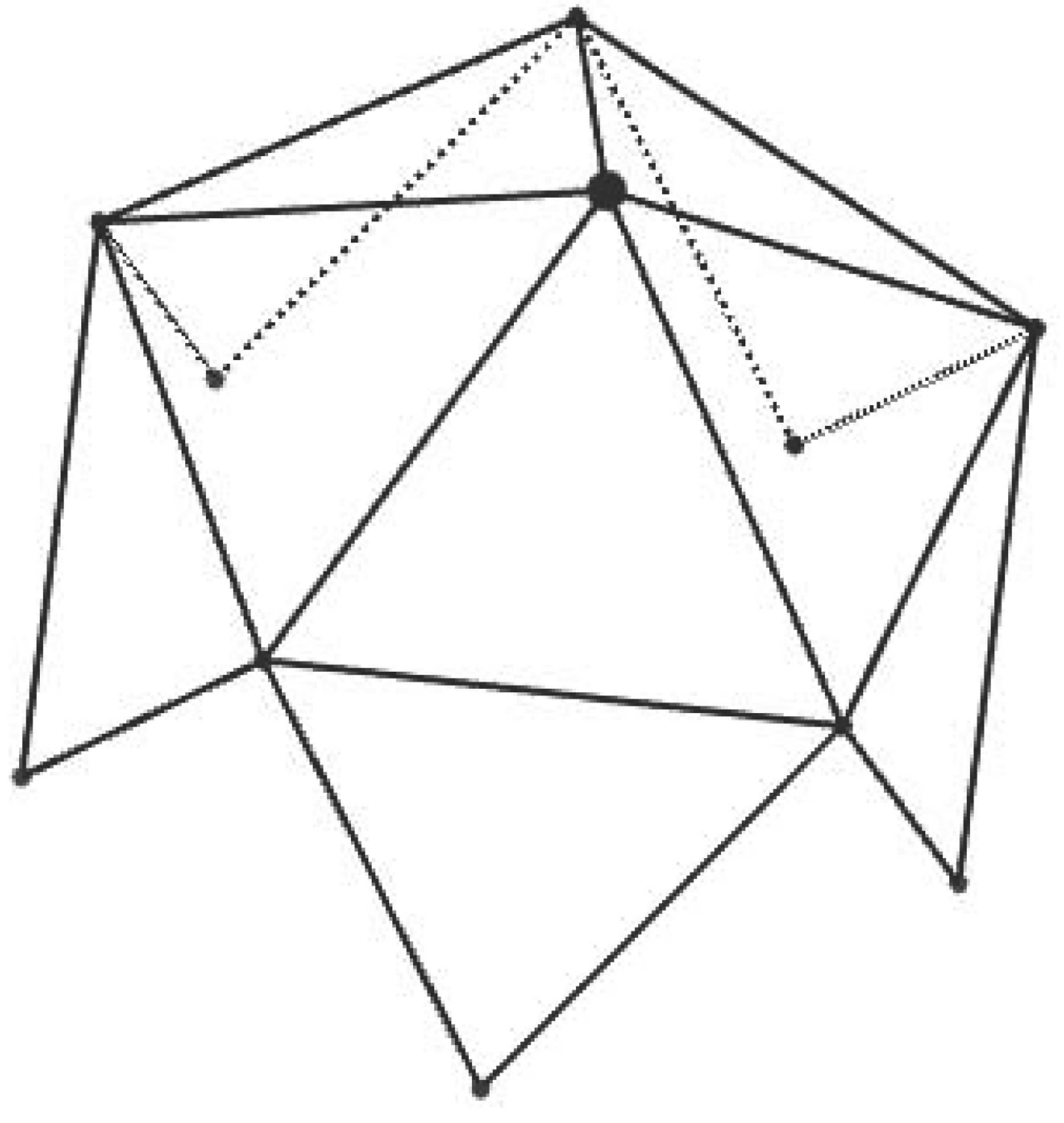}}
\end{minipage}}\par\medskip
\scalebox{.9}{\begin{minipage}{.32\linewidth}
\centering
\subfloat[]{\label{figure2d}\includegraphics[scale=.35]{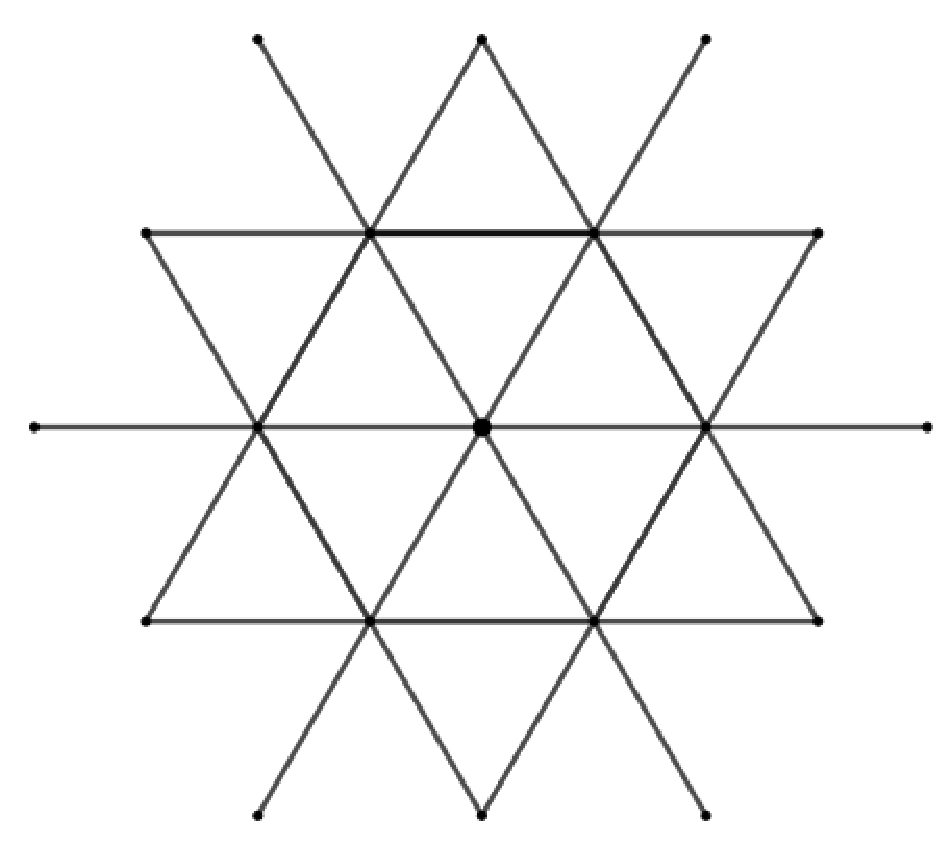}}
\end{minipage}%
\begin{minipage}{.35\linewidth}
\centering
\subfloat[]{\label{figure2e}\includegraphics[scale=.33]{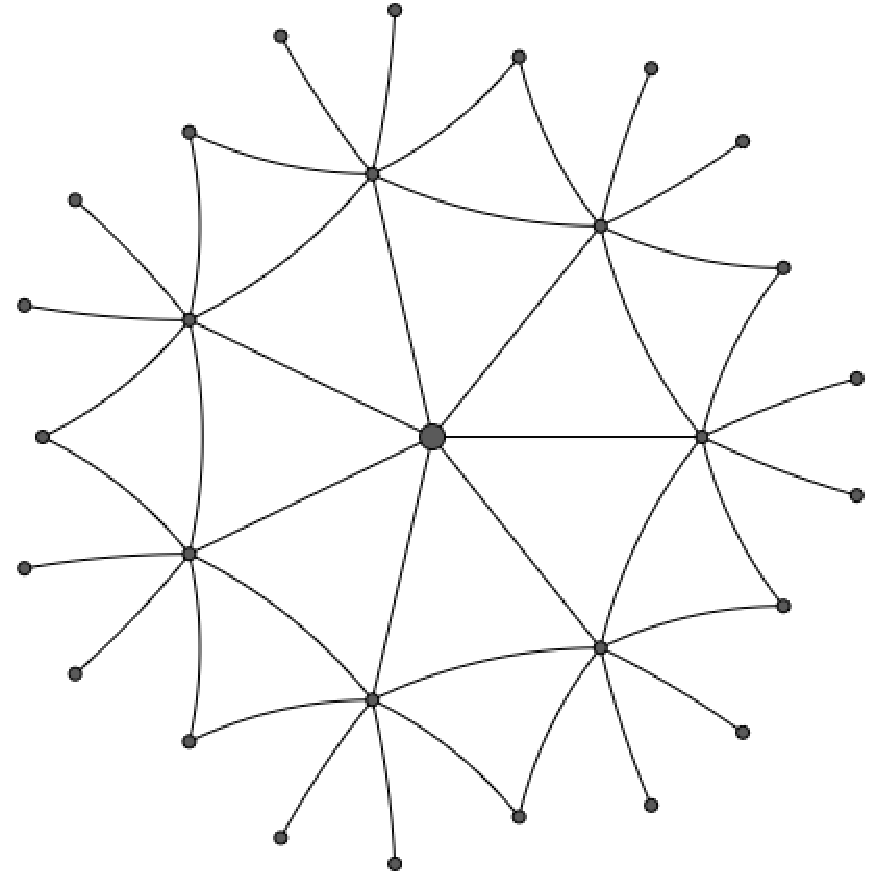}}
\end{minipage}%
\begin{minipage}{.32\linewidth}
\centering
\subfloat[]{\label{figure2f}\includegraphics[scale=.28]{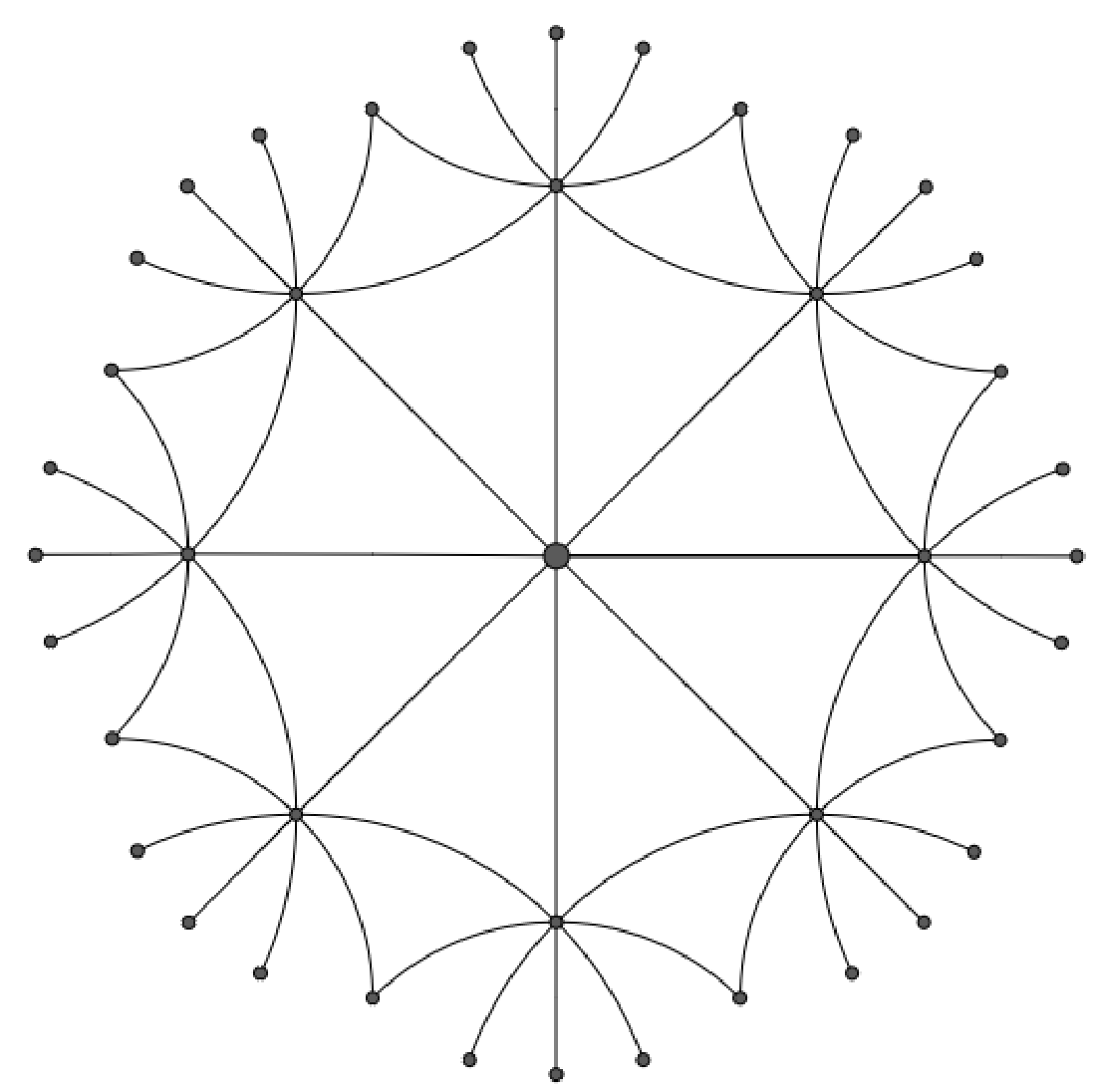}}
\end{minipage}}\par\medskip
\centering
\scalebox{.9}{\subfloat[]{\label{figure2g}\includegraphics[scale=.25]{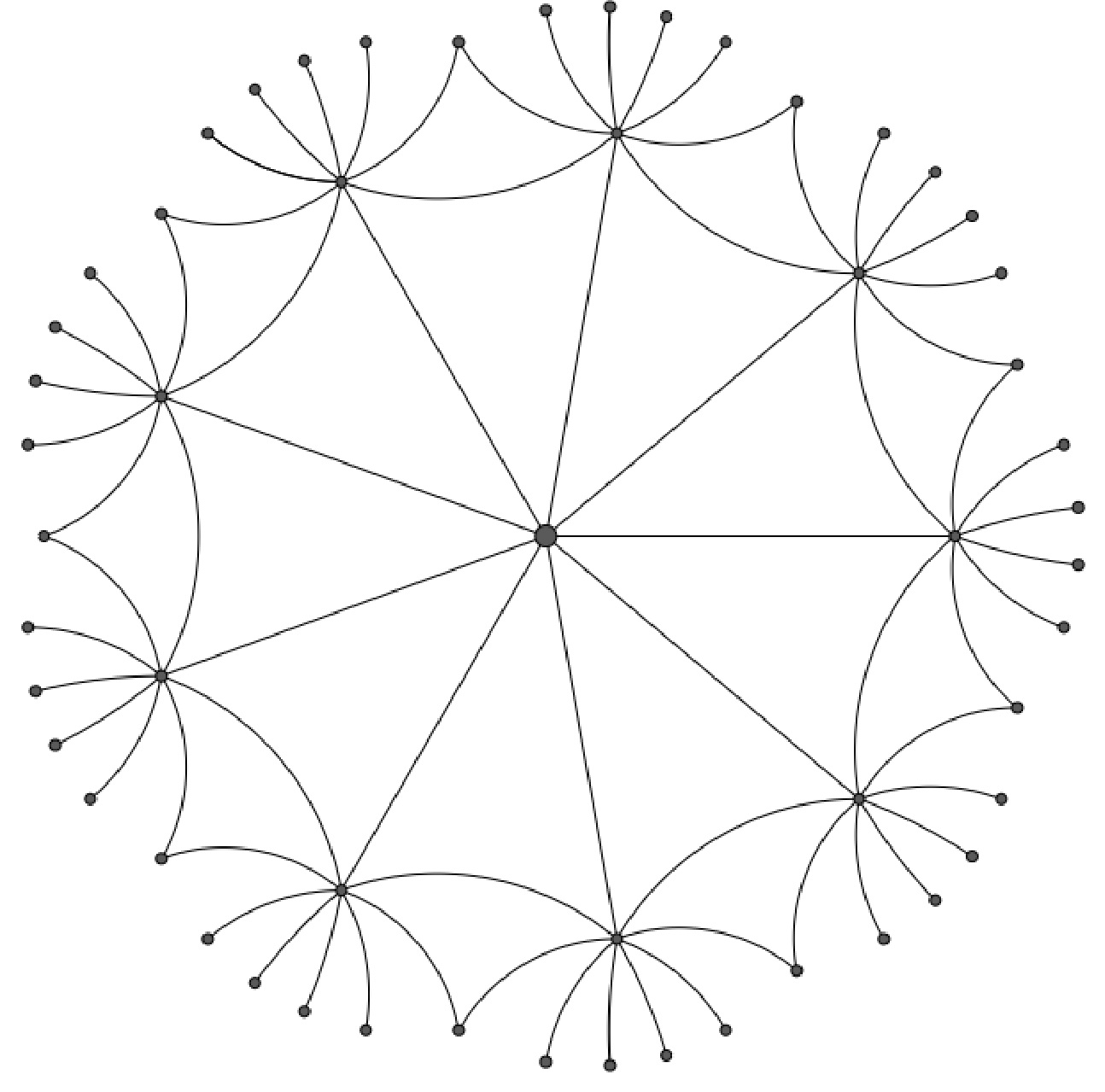}}}
\caption{The 2-step neighbourhoods of the considered graphs.}
\label{figure2}
\end{figure}

In Figure \ref{figure2}, we show the 2-neighbourhoods of a vertex of the corresponding graphs which are actually used in the discrete curvature computations.

The results are shown in Table \ref{curvaturetable}. In the second column of this table, the smooth curvatures of the associated ambient space forms are given. In the third column of the table, we give the discrete curvatures of the corresponding graphs computed by the program of Cushing et. al. \cite{CushingCalc2022}.

\begin{table}[H]
\centering
\begin{tabular}{ | c | c | c | }
\hline
\makecell{Type of the graph} & \makecell{The curvature of the \\ embedded space form} & \makecell{The discrete curvature\\ of the graph} \\ \hline
 \makecell{Order-3 regular triangular tiling\\ of the sphere} & 3.651 & 3 \\ \hline
 \makecell{Order-4 regular triangular tiling\\ of the sphere} & 2.467 & 3 \\ \hline
  \makecell{Order-5 regular triangular tiling\\ of the sphere} & 1.226 & 2.146 \\ \hline
  \makecell{Order-6 regular triangular tiling\\ of the plane} & 0 & 0 \\ \hline
 \makecell{Order-7 regular triangular tiling\\ of the hyperbolic plane} & -1.189 & -1.741 \\ \hline
  \makecell{Order-8 regular triangular tiling\\ of the hyperbolic plane} & -2.337 & -3.243 \\ \hline
 \makecell{Order-9 regular triangular tiling\\ of the hyperbolic plane} & -3.441 & -4.596 \\
 \hline
\end{tabular}\caption{} \label{curvaturetable}
\end{table}

The good news is that, for this collection of examples, there is at least complete coincidence of the signs of curvatures. Moreover, from the order-3 regular triangular tiling of the sphere downwards to the order-9 regular triangular tiling of the hyperbolic plane, the curvatures of the embedded space forms (second column of Table \ref{curvaturetable}) constitute a decreasing sequence. Remarkably, the discrete curvatures of the corresponding graphs (third column of Table \ref{curvaturetable}) also show the same pattern. The direct comparison in each row of Table \ref{curvaturetable} (the curvature of the embedded space form versus the discrete curvature of the corresponding graph) might not be very satisfactory, but that the discrete curvature somehow senses the smooth one seems rather reasonable. Such a correlation, if underpinned more thoroughly, would give more weight to the notion of the discrete curvature and any theoretical arguments, hints or heuristics in this direction would be very welcome.

\par


\end{document}